\documentclass[10pt]{amsart}
\usepackage{amsfonts,amsmath,amssymb}

\numberwithin{equation}{section}
\newtheorem{theo}[equation]{Theorem}
\newtheorem{pro}[equation]{Proposition}
\newtheorem{cor}[equation]{Corollary}
\newtheorem{lem}[equation]{Lemma}
\newtheorem{defi}[equation]{Definition}
\newtheorem{rem}[equation]{Remark}

\def\mg{\mathfrak {g}}
\def\k{\mathbb K}
\def\c{\mathbb C}
\def\z{\mathbb Z}

\def\t{\otimes}
\def\lra{\longrightarrow}
\def\H{\operatorname{H}}
\def\HL{\operatorname{HL}}

\def\B{\operatorname{B}}
\def\dim{\operatorname{dim}}

\def\m{\mathcal}
\def\d{\operatorname{d}}
\def\Hom{\operatorname{Hom}}
\def\Z{\mathbb {Z}}
\def\im{\operatorname{Im}}
\def\st{\stackrel}
\def\al{\alpha}
\def\be{\beta}
\def\ga{\gamma}
\title[Leibniz central extensions of Lie algebras]
{A cohomological characterization of Leibniz central extensions of Lie algebras}

\author[Hu]{Naihong Hu$^\star$}
\address{Department of Mathematics, East China Normal University,
Shanghai 200062, PR China} \email{nhhu@euler.math.ecnu.edu.cn \&
nhu@ictp.it (until August 31, 2006)}\author[Pei]{Yufeng Pei}
\address{Department of Mathematics, East China Normal University,
Shanghai 200062, PR China}\email{peiyufeng@gmail.com}
\author[Liu]{Dong Liu}
\address{Department of Mathematics, Jiaotong University,
Shanghai 200030, PR China}\email{dliu@sjtu.edu.cn}
\thanks{$^\star$N.H., Corresponding author. This work is
supported in part by the NNSF (Grants 10431040, 10671027), the
TRAPOYT and the FUDP from the MOE of China, the SRSTP from the
STCSM}

\begin{document}

\begin{abstract}
Mainly motivated by Pirashvili's spectral sequences on a Leibniz
algebra, a cohomological characterization of Leibniz central
extensions of Lie algebras is given based on Corollary 3.3 and
Theorem 3.5. In particular, as applications, we obtain the
cohomological version of Gao's main Theorem in \cite{Gao2} for
Kac-Moody algebras and answer a question in \cite{LH}.
\end{abstract}

\maketitle
\section{Introduction}
Leibniz algebras introduced by Loday (\cite{Lo}) are a
non-antisymmetric generalization of Lie algebras. There is a
(co)homology theory for these algebraic objects whose properties are
similar to those of the classical Chevalley-Eilenberg {co}homology
theory for Lie algebras. Since a Lie algebra is a Leibniz algebra,
it is interesting to study Leibniz (co)homology of Lie algebras
which may provide new invariants for Lie algebras. Lodder
(\cite{Lo1}) obtained the Godbillon-Vey invariants for foliations by
computing Leibniz cohomology of certain Lie algebras and mentioned
how a Leibniz algebra arises naturally from vertex (operator)
algebras. Recently, some interrelations with manifolds were
investigated, which could lead to possible applications of Leibniz
(co)homology in geometry (see \cite{HM, Lo1}, etc.).

Central extensions play a central role in the theory of Lie algebras
(see \cite{ACKP,B1,B,F0,F,Gao1,KPS,Ku,Li,Pa,W}, etc.). Viewing a Lie
algebra as a Leibniz algebra, it is natural to determine its Leibniz
central extensions, and compare the differences between Leibniz and
Lie central extensions. Loday-Pirashvili (\cite{LP}) have shown that
the Virasoro algebra is a universal central extension of the Witt
algebra in the category of Leibniz algebras as well.  It is
well-known that any Kac-Moody Lie algebra ${\frak g}(A)$ is
centrally closed (\cite{Gao1}). However, its universal central
extension in the category of Leibniz algebras is not centrally
closed for affine type as witnessed by Gao in (\cite{Gao2}). Other
cases of infinite dimensional Lie algebras were discussed in
\cite{LH, ZM}. Their approach involves technical and lengthy
computations.

The purpose of this note is to give a concise determination of the
universal central extension of a class of Lie algebras in the
category of Leibniz algebras. Our motivation comes from the
 long exact sequence  used by Pirashvili (\cite{Lo1, P}). Note that
the invariant symmetric bilinear forms and Leibniz central
extensions are connected by the Cartan-Koszul homomorphism (see
\cite{Ku}). In practice, it is easier to handle the invariant
symmetric  bilinear forms on a Lie algebra than to directly
determine its Leibniz central extensions. On the other hand, we
notice that determining Leibniz central extensions of a Lie algebra
is equivalent to treating dual space derivations of this Lie algebra
(cf. \cite{F}). We combine these observations to study the
(universal) Leibniz central extensions of Lie algebras. The current
method avoids complicated computations when applied to some classes
of infinite dimensional Lie algebras we are interested in.

The organization of this note is as follows. Section 2 introduces
some basic notions on Leibniz algebras. In Section 3,  we present a
short exact sequence as a variation of the Pirashvili's long exact
sequence. We then go on to describe the  Leibniz central extensions
of Lie algebras,  in terms of the invariant symmetric bilinear forms
and the Cartan-Koszul homomorphism (see Corollary 3.3),  and to
obtain
 the second Leibniz cohomology group with trivial coefficients by
dual space derivations (Theorem 3.5). Section 4 provides some
applications based on Corollary 3.3. As a consequence, a
cohomological version of Gao's result (\cite{Gao2}) in Kac-Moody
algebras case (Theorem 4.18) is given. In Section 5, we study the
Leibniz central extensions of the quadratic Lie algebras and
construct a counterexample to address a question posed in \cite{LH}.

\section{Prerequisites on Leibniz algebras}
\smallskip
\subsection{Leibniz algebra} Let $\k$ be an algebraically closed field with
$\operatorname{char}\k=0$.
\begin{defi}
 A Leibniz algebra is a $\k$-module $L$  with a bilinear
map $[-,-]:L \times L \longrightarrow L$  satisfying  the Leibniz
identity $[x, [y, z]]= [[x, y], z]-[[x, z], y]$, for $x, y, z\in L.
$
\end{defi}
The center of $L$ is defined as $\{z\in L\,|\,[z, L]=[L, z]=0\}.$
$L$ is called {\it perfect} if $[L,L]=L$. If, in addition,
$[x,x]=0$, $\forall\,x\in L$, the Leibniz identity is equivalent to
the Jacobi identity. In particular, Lie algebras are examples of
Leibniz algebras.

\begin{defi}
Let $L$ be a Leibniz algebra over $\k$.  $M$ is called a
representation of $L$ if $M$ is a $\k$-vector space equipped with
two actions $($left and right$)$ of $L$, i.e., $ [-,-]:L\times M\lra
M $ and $ [-,-]: M\times L\lra M $ satisfying
\begin{gather*}
(MLL)\qquad[m, [x, y]]= [[m, x], y]-[[m, y], x],
\\
(LML)\qquad[x, [m, y]]= [[x, m], y]-[[x, y], m],
\\
(LLM)\qquad[x, [y, m]]= [[x, y], m]-[[x, m], y],
\end{gather*}
for any $m\in M$ and $x,y\in L.$
\end{defi}

\subsection{Cohomology of Leibniz algebras}

Let $L$ be a Leibniz algebra over $\k$, and $M$ a representation of
$L$. Denote $C^n(L,M):=\operatorname{Hom}_{\k}( L^{\t n},M),\,n\geq
0$. The Loday coboundary map $ d^n: C^n(L,M)\rightarrow
C^{n+1}(L,M)$ is defined by
\begin{equation*}
\begin{split}
(d^nf)&(x_1,\cdots, x_{n+1})=[x_1,f(x_2,\cdots, x_{n+1})] \\
&+ \sum_{ i=2}^{n+1}(-1)^{i}[f( x_1,\cdots,\hat{x}_i,\cdots, x_{n+1}),x_i]\\
&+\sum_{1\le i<j\le n+1}(-1)^{j+1}f( x_1,\cdots, x_{i-1},
[x_i,x_j], x_{i+1}\cdots,\hat{x}_j,\cdots, x_{n+1}) .
\end{split}
\end{equation*}
Clearly, $d^{n+1}d^{n}=0,\ n\geq 0$. $(C^*(L,M),d)$ is a
well-defined cochain complex, whose cohomology is called the
cohomology of the Leibniz algebra $L$ with coefficients in the
representation $M$: $\HL^*(L,M):=\H^*((C^*(L,M),d))$. Similarly, we
have a chain complex $(C_*(L,M),d)$, whose homology is called the
homology of the Leibniz algebra $L$ with coefficients in the
representation $M$: $\HL_*(L,M):=\H_*((C_*(L,M),d))$.

\subsection{Cohomology of Lie algebras}
Let $\mathfrak{g}$ be a Lie algebra over $\k$, and $M$ a
$\mg$-module. Denote the Chevalley-Eilenberg cochain complex by
$$
(\Omega^*(\mg,M),\, \delta):=(\operatorname{Hom}( \wedge^*
\mg,M),\,\delta),
$$
where $\delta$ is the Chevalley-Eilenberg coboundary map defined by
\begin{equation*}
\begin{split}
(\delta^nf)&(x_1,\cdots, x_{n+1})
= \sum_{ i=1}^{n+1}(-1)^{i+1}x_i\cdot f( x_1,\cdots,\hat{x}_i,\cdots, x_{n+1})\\
&+\sum_{1\le i<j\le n+1}(-1)^{i+j}f([x_i,x_j], x_1,\cdots,
\hat{x}_{i} ,\cdots,\hat{x}_j,\cdots, x_{n+1}).
\end{split}
\end{equation*}
Then $\H^*(\mg,M):=\H^*((\Omega^*(\mg,M),\delta))$ is called the
cohomology of the Lie algebra $\mg$ with coefficients in the
$\mg$-module $M$.
\section{Leibniz cohomology and Leibniz central extensions of Lie algebras}

\subsection{Central extensions of Leibniz algebras}
A central extension of $L$ is a pair $(\hat{L},\pi)$, where
$\hat{L}$ is a Leibniz algebra and $\pi:\hat{L}\lra L$ is a
surjective homomorphism whose kernel lies in the center of
$\hat{L}$. The pair $(\hat{L},\pi)$ is a universal central extension
of $L$ if for every central extension $(\tilde{L},\tau)$ of $L$,
there is a unique homomorphism $\psi: \hat{L}\lra\tilde{L}$ for
which $\tau\circ\psi=\pi$. The following result is known.

\smallskip

\begin{pro}$($\cite{LP}$)$
There exists a one-to-one correspondence between the set of
equivalent classes of one-dimensional Leibniz central extensions of
$L$ by $\k$ and the second Leibniz cohomology group $\HL^2(L,\k)$.
\end{pro}

\subsection{Invariant symmetric bilinear forms}
Let $\mathfrak{g}$ be a Lie algebra over $\k$, and $M$ a
$\mg$-module. A symmetric bilinear form $\phi$ on $\mg$ is called
$\mg$-invariant if $\phi$ satisfies $\phi([x,y],z)=\phi(x,[y,z]),\
\forall\ x,y,z\in\mg$. Let $\B(\mg, \k)$ stand for the set of all
$\k$-valued symmetric $\mg$-invariant bilinear forms on Lie algebra
$\mg$.

The short exact sequence below is a variation of the first $5$-terms
of Pirashvili's long exact sequence (see \cite{P}).
\begin{pro} For any Lie algebra $\mg$, there is an exact sequence
$$
0\lra{\H^2(\mg,\k)}\st{f}{\lra}{\HL^2(\mg,\k)}\st{g}{\lra}{\B(\mg,
\k)}\st{h}{\lra}{\H^3(\mg,\k)},
$$
where
$f$ is the natural embedding map,
$g$ is defined by $ g(\alpha)(x,y)=\alpha(x,y)+\alpha(y,x)$,
$\forall\; \alpha\in\HL^2(\mg,\k)$, $x,\,y\in\mg$, and
$ h$ is the Cartan-Koszul map $($\cite{Ku}$)$ defined by $
h(\alpha)(x,y,z)=\alpha([x,y],z)$, $\forall\;\alpha\in\B(\mg,k)$,
$\,x,\, y,\, z\in\mg$.\hfill\qed
\end{pro}

\begin{cor} $\frac{\HL^2(\mg,\k)}{\H^2(\mg,\k)}=\ker(h)$. In
particular, $\HL^2(\mg,\k)=\H^2(\mg,\k)$ if and only if $\ker(h)=0$.
\end{cor}
\begin{rem}  Note that the natural embedding $
\Omega^*(\mg,M)\hookrightarrow C^*(\mg,M)$ induces a short exact
sequence in the category of cochain complexes of Lie algebra $\mg$:
$$
0\longrightarrow\Omega^*(\mg,M)\longrightarrow C^*(\mg,
M)\longrightarrow C^*_{rel}(\mg, M)[2]\longrightarrow 0,
$$ where
$C^*_{rel}(\mg,M)[2]:=\frac{C^*(\mg,M)} {\Omega^*(\mg,M)}$ is the
quotient cochain complex. Pirashvili's long exact sequence below
$($see \cite{Lo1}, \cite{P}$)$
$$
0\rightarrow \H^2(\mg,M)\rightarrow \HL^2(\mg,M)\rightarrow
\H^0_{rel}(\mg,M)\rightarrow \H^3(\mg,M)\rightarrow \HL^3(\mg,M)
\rightarrow \cdots$$ and isomorphisms $ \HL^i(\mg,M)= \H^i(\mg,M)$,
$i=0,\, 1$, are the main tools to compare the higher Lie and Leibniz
cohomology groups of a Lie algebra.

Let $M=\k$. We have
$$
0\rightarrow \H^2(\mg,\k)\rightarrow \HL^2(\mg,\k)\rightarrow
\H^0_{rel}(\mg,\k)\rightarrow \H^3(\mg,\k)\rightarrow \HL^3(\mg,\k)
\rightarrow \cdots. $$  As a consequence $($which is not clear for
us$)$ of the spectral sequences, Pirashvili claimed $($with no
proof, see \cite{P}\,$):$ $\H^0_{rel}(\mg,\k)\cong\B(\mg,\k)$, while
a direct elementary proof of Proposition 3.2 is given in the
Appendix.
\end{rem}

\subsection{Dual space derivations }
Let $\mg$ be a Lie algebra over $\k$, and $M$ a $\mg$-module. A
linear map $d: \mg\lra M$ is a derivation if $d([x,y])=x\cdot
d(y)-y\cdot d(x)$, for $\, x,\,y\in \mg$. The derivations of the
form $x\mapsto x\cdot m$ for some $m\in M$ are called inner
derivations.
 $\text{Der}(\mg,M)$ and $\text{Inn}(\mg,M)$ denote the spaces of
derivations and inner derivations, respectively. Clearly, $
\H^1(\mg,M)={\text{Der}(\mg,M)}/{\text{Inn}(\mg,M)}$ is the first
cohomology group of $\mg$ with coefficients in $M$.

\begin{theo}
$\HL^2(\mg,\k)=\H^1(\mg,\mg^*)=\operatorname{Der(\mg,\mg^*)}/\operatorname{Inn(\mg,\mg^*)}$,
for any Lie algebra $\mg$, where $\mg^*$ is the dual $\mg$-module.
\end{theo}
\begin{proof}
Define a map $\theta:\H^1(\mg,\mg^*)\lra\HL^2(\mg,\k)$ as
$\theta(\al)(x,y):=\al(y)(x)$ for any $\al\in\H^1(\mg,\mg^*)$ and
$x,\,y\in\mg$. Indeed, since
\begin{equation*}
\begin{split}
\theta(\al)(x,[y,z])&=\al([y,z])(x)=(y\cdot\al(z))(x)-(z\cdot\al(y))(x)\\
&=\al(z)([x,y])-\al(y)([x,z])=\theta(\al)([x,y],z)+\theta(\al)([z,x],y),
\end{split}
\end{equation*}
$\theta(\al)\in\HL^2(\mg,\k)$. $\theta$ is well-defined. If
$\theta(\al)=\bar 0$ in $\HL^2(\mg,\k)$, there exists an element
$\be$ in $C^1(\mg,\k)$ such that, for any $x,y\in\mg$, $
d^1(\be)(x,y)=\theta(\al)(x,y)$. That is, $-\be([x,y])=\al(y)(x)$.
So $\al$ is a $1$-coboundary in $\Omega^1(\mg,\mg^*)$, i.e.,
$\al=\bar 0\in\H^1(\mg,\mg^*).$ Hence, $\theta$ is injective. On the
other hand, for any $\be\in\HL^2(\mg,\k)$, we define a map
$\al\in\Hom_{\k}(\mg,\mg^*)$ by $\al(x)(y):=\be(y,x)$, for any
$x,y\in\mg$. Since
\begin{equation*}
\begin{split}
\al([x,y])(z)&=\be(z,[x,y])=\be([z,x],y)-\be([z,y],x)=\al(y)([z,x])-\al(x)([z,y]),\\
&=(x\cdot\al(x))(z)-(y\cdot\al(x))(z),\qquad\text{\it for any }\
z\in\mg,
\end{split}
\end{equation*}
we get $ \al([x,y])=x\cdot\al(x)-y\cdot\al(x) $ for any $x,y\in\mg$,
which means $\al\in\H^1(\mg,\mg^*)$. Hence, $\theta$ is surjective.
\end{proof}

Denote $\mbox{SDer}(\mg,\mg^*):=\{\phi\in\mbox{\rm
Der}(\mg,\mg^*)\;|\;\phi(x)(y)+\phi(y)(x)=0,\ \forall\
x,\,y\in\mg\}$. By Theorem 3.5 above and Proposition 1.3 (2) in
\cite{F0}, we have
\begin{cor}
$\HL^2(\mg,\k)/\H^2(\mg,\k)=\operatorname{Der(\mg,\mg^*)}/\operatorname{SDer(\mg,\mg^*)}\bigl(\subset
H^0_{rel}(\mg,\k)\bigr)$.
\end{cor}
\begin{cor}
If $\,\mg$ is a finite dimensional Lie algebra over $\k$ with a
nondegenerate invariant symmetric  bilinear form $\psi$, then $
\HL^2(\mg,\k)=\operatorname{Der(\mg,\mg)}/\operatorname{Inn(\mg,\mg)}.
$ If, in addition, $\mg$ is simple, then
$\HL^2(\mg,\k)=0$.\hfill\qed
\end{cor}

\section{Applications: Leibniz central extensions of some Lie algebras}

\subsection{Lie algebras of Virasoro type}

Let $\k[t,t^{-1}]$ be the Laurent polynomial algebra over $\k$ and
$\frac{d}{d t}$ be the differential operator on $\k[t,t^{-1}]$. Set
$L_n=-t^{n+1}\frac{d}{dt}$ and $I_n=t^n$ for $n\in\z$.
\begin{defi} \ The Witt algebra
$\mathcal{W}=\bigoplus_{n\in\mathbb {Z}}\mathbb {K}L_n $ is defined
by $[L_m,L_n]=(m-n)L_{m+n},\ m,\,n\in\mathbb {Z}$.
\end{defi}
It is well-known that
$\H^2(\mathcal{W},\k)=\H_2(\mathcal{W},\k)=\k\alpha$,
$\alpha(L_m,L_n)=\delta_{m+n},0\frac{m^3-m}{12}$.
\begin{defi} $($\cite{ACKP}$)$\ The Lie algebra
$\m{H}=\bigoplus_{n\in\z}\c L_n\bigoplus\bigoplus_{m\in \z}\c I_m$
is defined by $[L_m, L_n]=(m-n)L_{m+n}$,  $[I_m, I_n]=0$, $[L_m,
I_n]=-nI_{m+n}$, $m,\,n\in\mathbb {Z}$.
\end{defi}

In \cite{ACKP}, the authors proved that
$\dim\H^2(\m{H},\k)=\dim\H_2(\m{H},\k)=3.$ In fact, as the universal
central extension of $\m{H}$, the twisted Heisenberg-Virasoro
algebra $H_{Vir}$ has a basis $\{I_m, L_m, C_I, C_L, C_{LI}\mid
m\in\z\}$ satisfying the following relations:
\begin{gather*}
[I_m, I_n]=n\delta_{m+n, 0}C_{I},\\
[L_m, L_n]=(m-n)L_{m+n}+\delta_{m+n, 0}\frac{1}{12}(m^3-m)C_L,\\
[L_m, I_n]=-nI_{m+n}+\delta_{m+n, 0}(m^2-m)C_{LI},\\
[H_{Vir}, C_L]=[H_{Vir}, C_I]=[H_{Vir}, C_{LI}]=0.
\end{gather*}
\begin{defi}  Lie
algebra $\m{D}=\bigoplus_{m\in\z,n\in\z_{+}}\k t^{m}D^{n}$ of
differential operators is defined by
$[\,t^{m}D^{r},t^{n}D^{s}\,]=t^{m+n}((D+n)^r-(D+m)^s),\quad m,
\,n\in\z,\, r,\, s\in\z_+$, where $D=t\frac{d}{d t}$.
\end{defi}
By \cite{Li}, $\dim\H^2(\m{D},\k)=\dim\H_2(\m{D},\k)=1$, and the
universal central extension of $\m D$ is a Lie algebra
$\m{W}_{1+\infty}=\bigoplus_{m\in\z,n\in\z_{+}}\k
t^{m}D^{n}\bigoplus\k C$ with definition
$[t^{m}D^{r},t^{n}D^{s}]=t^{m+n}((D+n)^r-(D+m)^s)+\psi(t^mD^r,t^nD^s)C$,
$\psi(t^{m+r}D^r,t^{n+s}D^s)=\delta_{m+n,0}(-1)^rr!s!\binom{m+r}{r+s+1}$,
for $m,\,n\in\z,\, r,\,s\in\z_+$.
\begin{rem}
$\m{W}$ and $\m{H}$ are Lie subalgebras of $\m{D}$.
\end{rem}
\begin{pro} There is no non-trivial invariant symmetric  bilinear form
on Lie algebra $\mathcal{\mg}$, where $\mg=\m{W},\m{H}$, or $\m{D}$.
\end{pro}
\begin{proof} Assume that $f$ is an invariant symmetric bilinear form on
$\mathcal{\mg}$. Note that

$(1)\quad \mg=\m{W}$ is generated as Lie algebra by $L_{-2},\,L_3$
with $f(L_i,L_j)=0$ for $i,\,j\in\{-2,3\}$. So $f\equiv0$.

$(2)\quad \mg=\m{H}$  is generated as Lie algebra by $L_{-2},\,L_3$
and $I_1$ with $f(x, y)=0$ for $x, y\in\{L_{-2},L_3, I_1\}$. So
$f\equiv0$.

$(3)\quad \mg=\m{D}$  is generated as Lie algebra by $t,\,t^{-1}$
and $D^2$ (see \cite{Zh}). By a direct computation, we have
\begin{gather*}
t^{2}(\frac{d}{d t})^2=D^2-D,\\
t^{m}(\frac{d}{d t})^n=\frac{1}{n+1}\Bigl[t^m(\frac{d}{d
t})^{n+1},t\Bigr], \quad m\in \z,\, n\in\z_{+}.
\end{gather*}

By $(1)$ and $(2)$, $f(t^{\pm1},t^{\pm1})=f(t^{\pm
1},D)=f(D,D)=0.$
\begin{equation*}
\begin{split}
f\Bigl(t(\frac{d}{d t})^2,t^2(\frac{d}{d
t})^2\Bigr)&=f\Bigl(t(\frac{d}{d t})^2,-\frac{1}{3}[t,t^2(\frac{d}{d
t})^3]\Bigr){=-}\frac{1}{3}f\Bigl([t(\frac{d}{d
t})^2,t],t^2(\frac{d}{d
t})^3\Bigr)\\
&=-\frac{2}{3}f\Bigl(t\frac{d}{d t},t^2(\frac{d}{d
t})^3\Bigr)=-\frac{2}{3}f\Bigl(t\frac{d}{d
t},-\frac{1}{4}[t,t^2(\frac{d}{d t})^4]\Bigr)\\
&=\frac{1}{6}f\Bigl([t\frac{d}{d t},t],t^2(\frac{d}{d
t})^4\Bigr)
=\frac{1}{6}f\Bigl(t,t^2(\frac{d}{d t})^4\Bigr)\\
&=\frac{1}{6}f\Bigl(t,-\frac{1}{5}[t,t^2(\frac{d}{d
t})^5]\Bigr)=-\frac{1}{30}f\Bigl([t,t],t^2(\frac{d}{d t})^5\Bigr)\\
&=0.
\end{split}
\end{equation*}
Similarly, we have $f(D,t^2(\frac{d}{d t})^2)=0$. Then,
\begin{equation*}
\begin{split}
f(D^2,D^2)&=f\Bigl(D+t^{2}(\frac{d}{d t})^2,D+t^{2}(\frac{d}{d t})^2\Bigr)\\
&=f(D,D)+2f\Bigl(D,t^2(\frac{d}{d t})^2\Bigr)+f\Bigl(t^2(\frac{d}{d
t})^2,t^2(\frac{d}{d t})^2\Bigr)\\
&=2f\Bigl(D,t^2(\frac{d}{d
 t})^2\Bigr)+f\Bigl(t^2(\frac{d}{d t})^2,t^2(\frac{d}{d t})^2\Bigr)\\
 &=0.
\end{split}
\end{equation*}

In summary, $B(\mathfrak{g}, \mathbb K)=0$ for
$\frak{g}=\m{W},\,\m{H}$, or $\m{D}$.
\end{proof}
\begin{cor} $\;\HL^2(\mathcal{\mg},\mathbb {K})=\H^2(\mathcal{\mg},\mathbb {K})$,
for  $\mg=\m{W},\,\m{H}$, or
$\m{D}$.
\end{cor}
\begin{rem}
Corollary 4.6 for $\mg=\m{W}$ was obtained by \cite{LP}, for
$\mg=\m{D}$ by \cite{ LH}.
\end{rem}

\subsection{Lie algebras of Block type}

\begin{defi}$($\cite{DZ}$)$ Let $A$ be a torsion-free abelian group. $\phi: A\times
A\longrightarrow  \k $ is a non-degenerate, skew-symmetric,
$\Z$-bilinear function. Then the degenerate Block algebra
$\m{L}(A,\phi)=\bigoplus_{x\in A-\{0\}}\k e_x$ is defined by $
[e_x,e_y]=\phi(x,y)e_{x+y}$.
\end{defi}
\begin{defi} $($\cite{KPS}$)$
Let $A=\Z\times\Z$, $\phi({(m,n)},{(m_1,n_1)})=nm_1-mn_1$ be a
skew-symmetric bi-additive function. The Virasoro-like algebra $$
\mathcal{V}=\bigoplus_{(m,n)\in \z\times\z-\{(0,0)\}}\k e_{m,n}$$ is
defined by $ [e_{m,n},e_{m_1,n_1}]=(nm_1-mn_1)e_{m+n,m_1+n_1}$.
\end{defi}
\begin{defi} $($\cite{KPS}$)$
Let $A=\Z\times\Z$, $ \phi((m,n),(m_1,n_1))=q^{nm_1}-q^{mn_1}$ with
a fixed $q\in \k^*$ $($non-root of unity$)$, where $\phi$ is a
skew-symmetric bi-additive function. Then the $q$-analogue
Virasoro-like algebra $ \mathcal{V}_q=\bigoplus_{(m,n)\in
\z\times\z-\{(0,0)\}}\k e_{m,n} $ is defined by $
 [e_{m,n},e_{m_1,n_1}]=(q^{nm_1}-q^{mn_1})e_{m+n,m_1+n_1}$.
\end{defi}
\begin{lem}$($\cite{ZM}$)\  \B(\m{L}(A,\phi),\k)=\k\alpha$, where $
\alpha(e_x,e_y)=\delta_{x+y, 0}$.
\end{lem}
\begin{pro}
$\HL^2(\m{L}(A,\phi),\k)=\H^2(\m{L}(A,\phi),\k)$.
\end{pro}
\begin{proof} By Corollary 3.3 and Lemma 4.11, it suffices to prove that the
image of $\alpha$ under the Cartan-Koszul homomorphism $h$ is
non-zero. If $h(\alpha)=\bar 0\in \H^3(\mg,\k)$, there is a
$\psi\in\Lambda^2\mg$ such that
$h(\alpha)(e_x,e_y,e_z)=\alpha([e_x,e_y],e_z)=d(\psi)(e_x,e_y,e_z)
 =\psi(e_x,[e_y,e_z])+\psi(e_y,[e_z,e_x])+\psi(e_z,[e_x,e_y])
$, for $e_x, e_y, e_z \in\m{L}(A,\phi)$. Let $x+y+z=0$ and
$\phi(x,y)\neq 0$. Then $
\phi(x,y)\psi(e_x,e_{-x})+\phi(x,y)\psi(e_y,e_{-y})+\phi(x,y)\psi(e_z,e_{-z})=\phi(x,y)$,
that is, $\psi(e_x,e_{-x})+\psi(e_y,e_{-y})+\psi(e_z,e_{-z})=1$.  On
the other hand, since $-x-y-z=0$, the above identity holds for $-x,
-y, -z$. As $\psi$ is skew-symmetric, we have
$\psi(e_x,e_{-x})+\psi(e_y,e_{-y})+\psi(e_z,e_{-z})=-1$. It is
impossible. So $h(\alpha)\neq \bar 0$.
\end{proof}
 It follows from \cite{DZ} that
\begin{cor}
$$\HL^2(\m{L}(A,\phi),\k)=\{\,[\alpha_{\mu}]\mid \alpha_{\mu}(x,y)
=\delta_{x+y,0}\mu(x),\forall\
x,\,y\in A,\ \mu\in\operatorname{Hom}_{\k}(A,\k)\}.$$
\end{cor}
\begin{rem}
With a different method, Corollary 4.13 was given in \cite{ZM}. For
the $($$q$-analogue$)$ Virasoro-like algebras, the same result was
given in \cite{LH}.
\end{rem}

\subsection{Kac-Moody algebras}

\begin{lem}
For a Lie algebra $\mg$ with $\dim \B(\mg,\k)\le 1$, in both cases
below$:$
\begin{itemize}
\item if \ $\B(\mg,\k)=0;$ or
\item if  \ $\B(\mg,\k)=\k\phi$ $(\phi\ne0)$, and there exists a subalgebra $\mathfrak a\cong
\mathfrak{sl}(2,\k)$ such that $\phi|_{\mathfrak{a}}\ne 0$,
\end{itemize}
then $\HL^2(\mg,\k)=\H^2(\mg,\k).$
\end{lem}
\begin{proof} The first case is clear by Corollary 3.3. Now for nonzero $\phi\in\B(\mg,\k)$, if there is a subalgebra $\mathfrak{a}$ as a
$\mathfrak{sl}(2,\k)$-copy such that $\phi|_{\mathfrak{a}}\ne 0$,
then $\phi|_{\mathfrak{a}}$ is a nonzero scalar multiple of the
Killing form on $\mathfrak{a}$. Let $h$ be the Cartan-Koszul
homomorphism: $h(\phi)(a,b,c)=\phi([a,b],c)$ for any $a, b,
c\in\mg$. If $h(\phi)=\bar 0$, then there is a
$\theta\in\Lambda^2\mg$ such that $\delta^2(\theta)=h(\phi)$, i.e.,
$\theta(a,[b,c])+\theta(b,[c,a])+\theta(c,[a,b])=\phi([a,b],c)$, for
$ a,\,b,\,c\in\mg$. Take $x,\,y,\,h\in\mathfrak{a}$ satisfying
$[x,y]=h,\;[h,y]=-2y,\; [h,x]=2x.$ Then
\begin{gather*}
\theta(x,[y,h])+\theta(y,[h,x])+\theta(h,[x,y])=\phi([x,y],h),\\
2\theta(x,y)+2\theta(y,x)+\theta(h,h)=\phi(h,h).
\end{gather*}
Therefore, $\phi(h,h)=0$. It contradicts the property of the Killing
form.
\end{proof}
\begin{lem} $\HL^2(\mg,\k)=\operatorname{Hom}(\HL_2(\mg,\k),\k)$, for
any perfect Leibniz algebra $\mg$.
\end{lem}
\begin{proof} Following the universal coefficient theorem for Leibniz
algebras in \cite{CFV}, one has the following short exact sequence
$$
0\longrightarrow\operatorname{Ext}(\HL_1(\mg,\k),\k)\longrightarrow
\HL^2(\mg,\k)\longrightarrow\operatorname{Hom}(\HL_2(\mg,\k),\k)\longrightarrow
0,
$$
where $\operatorname{Ext}(\HL_1(\mg,\k),\k)$ denotes the abelian
extension of $\HL_1(\mg,\k)$ by $\k$. It is clear that
$\operatorname{Ext}(\HL_1(\mg,\k),\k)=0$ since
$\HL_1(\mg,\k)=\mg/[\mg,\mg]=0$.
\end{proof}

Let $\mg$ be a perfect Lie algebra $([\mg,\mg]=\mg)$ over $\k$.
Then, by \cite{G}, there are a universal central extension
$\pi:\hat{\mg}\twoheadrightarrow\mg$ in the category of Leibniz
algebras, and a universal central extension
$\tilde{\pi}:\tilde{\mg}\twoheadrightarrow\mg$ in the category of
Lie algebras.
\begin{lem}
$\HL_2(\tilde{\mg},\k)=\ker \{\HL_2(\mg,\k)\twoheadrightarrow\H_2(\mg,\k)\}.$
\end{lem}
\begin{proof}
See 4.6 in \cite{LP} or Corollary 2.7 in \cite{G}.
\end{proof}

Let $R$ be a unital, commutative and associative algebra over $\k$.
The $R$-module of K\"ahler differentials $\Omega_{R|\k}^{1}$ is
generated by $\k$-linear symbols $\d{a}$ for $a\in R$ with the
relation $\d(ab)=a\d b+b\d a$, for any $a,\, b\in R$. In particular,
if $R=\k[t,t^{-1}]$, then $\Omega_{R|\k}^{1}=\bigoplus_{m\in\z}\k
t^m \d{t}$. For any positive integer $r$, let
$$
\Omega_{R|\k}^{1}(r)=\bigoplus_{i\in\z}\k t^{ir-1} \d{t},\quad
\Omega_{R|\k}'^{1}(r)=\bigoplus_{i\in\z-\{0\}}\k t^{ir-1} \d{t}
$$ be two
$\k$-subspaces of $\Omega_{R|\k}^{1}$. Then $\Omega_{R|\k}'^{1}(r)$
is a subspace of $\Omega_{R|\k}^{1}(r)$ of codimension $1$.

\begin{theo} Let $A=(a_{ij})_{1\le i,j\le n}$ denote an $n\times n$-matrix of
rank $\ell$ with entries in  $\k$. Using the notations in \cite{K},
denote by $\mg(A)$ the Kac-Moody Lie algebra associated to $A$. Let
$\mg'(A)$ be the derived algebra of $\mg(A)$, $\mathfrak{c}$ the
center of $\mg'(A)$, and $\bar{\mg}(A)=\mg'(A)/\mathfrak{c}$. Thus,

$(1)$\ If $a_{ii}\ne 0,\ 1\le i\le n$, then
$\dim\HL^2(\mg(A),\k)=(n-\ell)^{2}$;

$(2)$\ If $A$ is an indecomposable generalized Cartan matrix of
affine $X_n^{(r)}$ type, then
\begin{gather*}
\HL^2(\bar{\mg}(A),\k)=(\Omega_{R|\k}^{1}(r))^*,\qquad
\HL^2({\mg}'(A),\k)=(\Omega_{R|\k}'^{1}(r))^*.
\end{gather*}

$(3)$\ If $A$ is an indecomposable generalized Cartan matrix of
non-affine type, then
\begin{gather*}
\HL^2(\bar{\mg}(A),\k)=\mathfrak{c^*},\qquad \HL^2({\mg}'(A),\k)=0.
\end{gather*}
\end{theo}
\begin{proof}
(1)\ Theorem 3.2 \cite{F} tells us that $
\dim\H^1(\mg(A),(\mg(A))^*)=(n-\ell)^2$ under the assumption
$a_{ii}\neq0$ for $1\le i\le n$. Then Theorem 3.5 gives the result.

(2)\ Let $A$ be an indecomposable generalized Cartan matrix of
affine $X_n^{(r)}$ type. By Gabber-Kac's radical Theorem (Theorem
9.11 and remarks in \cite{K}, pp. 159),  $\mg'(A)$ and
$\bar{\mg}(A)$ can be presented in term of $3n$ generators
$f_i,\,h_i,\,e_i$ ($1\le i\le n$) and the Chevalley-Serre relations.
By Theorem 3.17 in \cite{Gao2}\footnote{Note that author in
\cite{Gao2} used the different notations.}, which says
$\HL_2(\bar{\mg}(A),\k)=\Omega_{R|\k}^{1}(r)$, together with Lemma
4.16, we get $\HL^2(\bar{\mg}(A),\k)=(\Omega_{R|\k}^{1}(r))^*$.

Because $\mg'(A)$ is the universal covering of $\bar{\mg}(A)$ in the
category of Lie algebras (see \cite{W}), we have the following exact
sequence:
$$
0\lra\H_2(\bar{\mg}(A),\k)\lra\mg'(A)\lra\bar{\mg}(A)\lra 0,
$$
where
$\H_2(\bar{\mg}(A),\k)=\Omega_{R|\k}^{1}(r)/\Omega_{R|\k}'^{1}(r)=t^{-1}\d
t$. Lemma 4.17 then yields
$$ \HL_2({\mg'}(A),\k)=\ker \{\HL_2(\bar{\mg}(A),\k)\twoheadrightarrow\H_2(\bar{\mg}(A),\k)\}=\Omega_{R|\k}'^{1}(r).$$
So Lemma 4.16 gives the second result.

(3)\  Let $A$ be an indecomposable generalized Cartan matrix of
non-affine type. Recall Berman's Theorem 3.1 \cite{B}, which says
that $\bar{\mg}(A)$ possesses a non-degenerate invariant symmetric
bilinear form $\phi$ if and only if $A$ is symmetrizable. Moreover,
such forms on $\bar{\mg}(A)$ are unique up to scalars.

If $A$ is symmetrizable, then Berman's result above means
$\dim\B(\bar{\mg}(A),\k)\ge1$. On the other hand, for any
$0\ne\psi\in\B(\bar{\mg}(A),\k)$, $\psi$ is necessarily
non-degenerate because $\bar{\mg}(A)$ is simple (by Theorem 4.3
\cite{K}, Exercise 1.4 \cite{K}). So Berman's result above insures
$\B(\bar{\mg}(A),\k)=\k\phi$. The non-degeneracy of $\phi$ implies
the nontrivial property of its restriction to a
$\mathfrak{sl}(2,\k)$-copy (see \cite{B1}). Lemma 4.15 shows $
\HL^2(\bar{\mg}(A),\k)=\H^2(\bar{\mg}(A),\k)=\mathfrak{c^*}$ (the
2nd ``$=$" was proved by Theorem 2.3 \cite{B}).

If A is non-symmetrizable, there exists no non-degenerate invariant
symmetric bilinear form on $\bar{\mg}(A)$. Since $\bar{\mg}(A)$ is
simple (by Theorem 4.3 \cite{K}, Exercise 1.4 \cite{K}), there is no
invariant symmetric bilinear form on $\bar{\mg}(A)$, i.e.,
$\B(\bar{\mg}(A), \k)=0$. Then Corollary 3.3 or Lemma 4.15 gives
$\HL^2(\bar{\mg}(A),\k)=\H^2(\bar{\mg}(A),\k)=\mathfrak{c^*}$ (the
2nd ``$=$" was proved by Theorem 2.3 \cite{B}).

Furthermore, since $\mg'(A)$ is the universal Leibniz central
extension of $\bar{\mg}(A)=\mg'(A)/\mathfrak{c}$. Consequently,
$\HL^2({\mg}'(A),\k)=0.$
\end{proof}

\begin{rem}
Theorem 4.18 $(2)$, $(3)$ give a criterion to distinguish between
affine and non-affine Kac-Moody algebras by means of the vanishing
property of the second Leibniz cohomology groups of
$\mathfrak{g}'(A)$ with trivial coefficients, where the homological
versions of Theorem 4.18 $(2)$ and the second statement of $(3)$
were due to Gao $($\cite{Gao2}$)$. However, strictly speaking, owing
to Gabber-Kac's Theorem $($see \cite{K}$)$, the definition of
$\mathfrak{g}'(A)$ for non-symmetrizable cases adopted by Gao is
different from ours used here, since our $\mathfrak{g}'(A)$ in Kac's
notation is the quotient of the former. So in this sense, we get the
same result for the quotient object.
\end{rem}

\section{Quadratic Leibniz algebras and their central extensions}

\subsection{Quadratic Leibniz algebra}
\begin{defi}
$(\mg,\phi)$ is called a quadratic Leibniz algebra if $\phi$ is a
symmetric invariant bilinear form on Leibniz algebra $\mg$.
\end{defi}

\begin{lem}
Let $\mg$ be a Lie algebra. If $(\mg,\phi)$ is a quadratic Leibniz
algebra and $d$ is a derivation of $\mg$, then $f(x,y):=\phi(x,dy)$
is a Leibniz $2$-cocycle on $\mg$.
\end{lem}
\begin{proof} Since
$f(x,[y,z])=\phi(x,d[y,z])=f([x,y],z)-f([x,z],y)$ for $x,\, y,\,
z\in\mg$, $f$ is a Leibniz $2$-cocycle on $\mg$.
\end{proof}
\begin{cor} $($\cite{K}$)$
If d is a skew-derivation, $\operatorname{i.e.},\
\phi(dx,y)+\phi(x,dy)=0$, then $f(x,y):=\phi(x,dy)$ is a Lie
$2$-cocycle on $\mg$.
\end{cor}

\subsection{A negative answer to a question in \cite{LH} }
A question in \cite{LH} is posed: If each Leibniz $2$-cocycle on an
infinite-dimensional Lie algebra is also a Lie $2$-cocycle.
\begin{pro} Assume that $\mg$ is a finite dimensional
simple Lie algebra over $\k$. Construct the Lie algebra
$\mg\t\k((t))$ with bracket: $ [x\t r,y\t s]'=[x,y]\t rs,\
x,\,y\in\mg,\, r,\, s\in\k((t))$. Then
 $\mg\t\k((t))$ is an infinite-dimensional simple Lie algebra over $\k$ and $$
\H^2(\mg\t\k((t)),\k)\subsetneq\HL^2(\mg\t\k((t)),\k).
$$
\end{pro}
\begin{proof}
Define $\phi(x\t r,y\t s)=(x,y)\operatorname{Res}(rs)$, where
$\operatorname{Res}$ is a linear function on $\k((t))$ and takes the
coefficient of $t^{-1}$ for every series, $(\,,)$ is the Killing
form on $\mg$. Then $\phi$ is an invariant symmetric  bilinear form
on $\mg\t\k((t))$ and $ t^k\frac{d}{d t},k\in\z-{\{0\}}$ is a
derivation of $\k((t))$. By Lemma 5.2, we get a non-trivial Leibniz
$2$-cocycle of the Lie algebra $\mg\t\k((t))$:
$$
f(x\t \sum_{m\geq N}a_mt^m,\,y\t \sum_{n\geq
N}b_nt^n)=(x,y)\sum_{m,n\geq N} na_mb_n\delta_{m+n+k,0},\eqno(5.1)$$
where $f$ is well-defined since the summation in (5.1) is finite,
and $f$ is not skew symmetric, that is, $
\H^2(\mg\t\k((t)),\k)\subsetneq\HL^2(\mg\t\k((t)),\k)$.
\end{proof}
\begin{rem}
$\mg\t\k((t))$ is a simple Lie algebra $($see \cite{ZM0}$)$, and $
\dim\H^2(\mg\t\k((t)),\k)$ $=1$, but
$\dim\HL^2(\mg\t\k((t)),\k)=\infty$, which then leads to a negative
answer to the above question \cite{LH}.
\end{rem}

\section{Appendix: \ An elementary proof of the exact sequence}

\noindent {\it Proof of Proposition 3.2.} \ We shall prove the
following
\begin{enumerate}
\item The map $g$ is well-defined, that is, $
g(\alpha)=g(\alpha+d^1\beta)\in\B(\mg,\k)$ for any
$\alpha\in\HL^2(\mg,\k)$ and $\beta\in \HL^1(\mg,\k)$ since
$g(d^1\beta)\equiv 0$ by definition. Clearly, $g(\alpha)$ is
symmetric. Since, for any $x,y,z\in\mg$,
\begin{equation*}
\begin{split}
g(\alpha)([x,y],z)&=\alpha([x,y],z)+\alpha(z,[x,y])\\
&=\alpha(x,[y,z])+\alpha([x,z],y)+\alpha([z,x],y])-\alpha([z,y],x)\\
&=\alpha(x,[y,z])+\alpha([y,z],x)\\
&=g(\alpha)(x,[y,z]),
\end{split}
\end{equation*}
$g(\alpha)$ is $\mg$-invariant.

\item $\im(f)=\ker(g)$. For any $\alpha\in\H^2(\mg,\k)$ and any $x,y\in\mg$, we have
$$g(f(\alpha))(x,y)=g(\alpha)(x,y)=\alpha(x,y)+\alpha(y,x)=0,$$
since $\alpha$ is skew symmetric. Then $\im(f)\subseteq\ker(g)$. On
the other hand, for any $\beta\in\ker(g)$, we have $ g(\beta)(x,y)=0
$, which implies that $\alpha(x,y)+\alpha(y,x)=0$ for any
$x,y\in\mg$. Then $\alpha$ is a skew symmetric Leibniz $2$-cocycle,
namely, a Lie $2$-cocycle in $\H^2(\mg,k)$. Hence, we have
$\ker(g)\subseteq\im(f).$

\item $\im(g)=\ker(h)$. Now, $\forall\;\alpha\in\HL^2(\mg,\k)$, define a bilinear form $\al'$ on $\mg$ by
$$\al'(x,y)=\al(y,x),$$
for any $x,y\in\mg$. Denote $\be:=\frac{\al+\al'}{2}$ and
$\ga:=\frac{\al-\al'}{2}$. Then $\be$ is symmetric and $\ga$ is skew
symmetric. Also, we have $\al=\be+\ga$.

Now, for any $x,y,z\in\mg$,
\begin{equation*}
\begin{split}
(h(g(\al)))(x,y,z)&=g(\al)([x,y],z)
=\al([x,y],z)+\al(z,[x,y])\\
&=\al([x,y],z)+\al([z,x],y)+\al([y,z],x)\\
&=\al(z,[x,y])+\al(y,[z,x])+\al(x,[y,z])\\
&\quad -\al([x,y],z)-\al([z,x],y)-\al([y,z],x)\\
&=2\ga(z,[x,y])+2\ga(y,[z,x])+2\ga(x,[y,z])\\
&=\delta^2(2\ga)(x,y,z),
\end{split}
\end{equation*}
which implies that $h(g(\al))$ is a $3$-coboundary in
$\Omega^3(\mg,\k)$. Then $\im(g)\subseteq\ker(h).$ Conversely, for
any $\omega\in \ker(h)$, there exists an element $\eta\in
\Omega^2(\mg,\k)$ such that $\delta^2(\eta)=h(\omega).$ Namely, for
any $x,y,z\in \mg$, there holds
$$
\eta(z,[x,y])+\eta(y,[z,x])+\eta(x,[y,z])=\omega([x,y],z).$$  Then
$$
\omega([x,y],z)+\eta([x,y],z)=\omega(x,[y,z])+\eta(x,[y,z])+\omega([x,z],y)+\eta([x,z],y),
$$
since $\omega(x,[y,z])+\omega([x,z],y)=0$. Let
$\rho=\frac{1}{2}(\omega+\eta)$. Then $ g(\rho)=\omega$ and
$$
\rho([x,y],z)=\rho(x,[y,z])+\rho([x,z],y),
$$ which means that $\rho\in\HL^2(\mg,\k)$.
Hence $\ker(h)\subseteq\im(g).$

We complete the proof.\hfill\qed
\end{enumerate}

\newpage \centerline{\bf ACKNOWLEDGMENT}

\vskip10pt The authors are indebted to the referee for the helpful
comments.  NH would like to express his thanks to the ICTP for its
support and hospitality when he visited the ICTP Mathematics Group,
Trieste, Italy, from March 1st to August 31st, 2006. YP is grateful
to Hagiwara and Papi for sending him their papers (\cite{HM},
\cite{Pa}).

\bibliographystyle{amsalpha}

\end{document}